\documentclass[11pt]{article}

\usepackage[T2A]{fontenc}
\usepackage[cp1251]{inputenc}
\usepackage{amsfonts}
\usepackage{eufrak}
\usepackage{amssymb}
\usepackage{amsmath}


\begin{document}

 \newtheorem{theorem}{Theorem}
\newtheorem{lemma}{Lemma}
\newtheorem{corollary}{Corollary}
\newtheorem{definition}{Определение}
\newtheorem{remark}{Remark}

\noindent The final version of the article will be published in the journal

\noindent Aequationes mathematicae

\vskip 1 cm

\noindent{\bf On a generalisation of the Skitovich--Darmois theorem}

\noindent{\bf for several linear forms on Abelian groups}

\bigskip

\noindent{Gennadiy Feldman}

\bigskip

\noindent{\bf Abstract.}  A.M. Kagan introduced a class of distributions $\mathcal{D}_{m, k}$ in $\mathbb{R}^m$ and proved that if the joint distribution of $m$ linear forms of $n$ independent random variables belongs to the class $\mathcal{D}_{m, m-1}$, then the random variables are Gaussian. A.M. Kagan's theorem implies, in particular, the well-known  Skitovich--Darmois theorem, where the Gaussian distribution on the real line is characterized by   independence of two linear forms of $n$ independent random variables.  In the note we describe a wide class of locally compact Abelian groups  where  A.M. Kagan's theorem is valid.

\bigskip

\noindent {\bf Mathematical Subject Classification:} 43A25, 43A35, 60B15, 62E10

\bigskip

\noindent {\bf Keywords:} locally compact Abelian group, Gaussian distribution, linear forms

 \bigskip

\noindent{\textbf{1. Introduction}}

\bigskip

\noindent One of the most famous characterization theorems in mathematical statistics is the following statement. It was proved independently by  V. P. Skitovich and G. Darmois.

\bigskip

\noindent\textbf{The Skitovich--Darmois theorem} (\cite[Chapter 3]{Kag-Lin-Rao}).
\textit{Let $\xi_i$, $i=1, 2,\dots, n$, $n\geq 2$,
  be independent random variables, $\alpha_i, \beta_i$  be nonzero real numbers.
If the linear forms $L_1=\alpha_1\xi_1+\cdots+\alpha_n\xi_n$ and
$L_2=\beta_1\xi_1+\cdots+\beta_n\xi_n$ are independent, then all the random variables   $\xi_i$ are Gaussian.}

\bigskip

  Consider  $m$ linear forms of $n$ independent random variables
$L_j=\alpha_{j1}\xi_1+\cdots+\alpha_{jn}\xi_n$, $j=1, 2, \dots, m$. In  \cite{Kag} A.M. Kagan introduced a class  of distributions   $\mathcal{D}_{m, k}$ in $\mathbb{R}^m$ and proved that if all the coefficients $\alpha_{ji}$ are nonzero and the joint distribution of the  linear forms $L_j$   belongs to the class $\mathcal{D}_{m, m-1}$, then the random variables $\xi_i$ are Gaussian.  In particular, if the linear forms are independent, then their joint distribution   belongs to the class $\mathcal{D}_{2,  1}$,  and hence, the Skitovich--Darmois theorem follows from  A.M. Kagan's theorem.  The aim of this note is to describe a wide class of locally compact Abelian groups where  A.M. Kagan's theorem is valid. In do doing coefficients of the linear forms are continuous endomorphisms of a group.  We  remark that a number of papers is devoted to group analogues of the Skitovich--Darmois theorem, see e.g. \cite{Fe1990}--\cite{Fe-SD-2003}, \cite{Fe2009}--\cite{MyFe2011}, and also \cite[\S10--15]{Fe1}.

We will use in the article standard results of abstract harmonic analysis (see \cite{Hewitt-Ross}). Recall some definitions and agree on notation.

Let $X$ be a locally compact Abelian group,
   $Y$ be its character group,
$(x,y)$ be the value of a character  $y\in Y$ at an element  $x\in X$.    If   $H$ is a subgroup of the group $Y$, denote its annihilator by   $A(X,H)=\{x \in X:
(x,y)=1$ for all  $y \in H\}$.
    Let $X_1$ and $X_2$
be locally compact Abelian groups, and  $Y_1$ and $Y_2$ be their character groups respectively.  For any continuous homomorphism   $f : X_1 \mapsto X_2$ define the adjoint homomorphism    $\tilde{f} :Y_2 \mapsto Y_1$ by the formula
$(x_1,\tilde{f}y_2)=(f x_1,y_2)$ for all $x_1 \in X_1$, $y_2
\in Y_2$. If $G$ is a subgroup of $X$, denote  its closure by $\overline G$.  Denote by $\mathbb{T}$ the circle group (the one dimensional torus), i.e.
${\mathbb{T}=\{z\in\mathbb{C}:\ |z|=1\}}$. Denote by $I$ the identity automorphism of a group.

Let $\psi(y)$ be a function on the group    $Y$,   and let $h \in
Y$. Denote by   $\Delta_h$   the finite difference operator
$$
\Delta_h \psi(y)=\psi(y+h)-\psi(y).
$$
A function $\psi(y)$ on $Y$ is called
a    polynomial  if
\begin{equation}\label{21}
 \Delta_{h}^{n+1}\psi(y)=0, \ y,h \in Y,
\end{equation} for some nonnegative integer  $n$.
The minimal $n$ for which (\ref{21}) holds is called   the
degree  of the polynomial $\psi(y)$.

Let $\xi$ be a random variable with values in the group $X$. Denote by $\mu_\xi$ its distribution and by $\hat \mu(y)=\hat \mu_\xi(y)$ its characteristic function (Fourier transform)
 $$\hat \mu(y)=\hat \mu_\xi(y)={  \mathbb{E}}[(\xi, y)]= \int_X (x, y) d\mu_\xi(x), \quad y\in Y,$$
  where ${  \mathbb{E}}$ is the    mathematical expectation of a complex valued random variable. Denote by ${\rm M}^1(X)$   the convolution semigroup of probability distributions on the group   $X$.
For $\mu \in{\rm M}^1(X)$ define a distribution   $\bar \mu \in
{\rm M}^1(X)$ by the formula $\bar \mu(B) = \mu(-B)$ for all Borel sets
   $B$. We note that
$\hat{\bar{\mu}}(y)=\overline{\hat\mu(y)}$.

A distribution  $\gamma$  on the group $X$ is called Gaussian
(\cite[Chapter IV]{Parthasarathy})
if its characteristic function can be represented in  the form
$$
\hat\gamma(y)=(x,y)\exp\{-\varphi (y)\}, \quad y\in Y,
$$
where $x \in X$, and $\varphi(y)$ is a continuous non-negative function
on the group $Y$
 satisfying the equation
$$
    \varphi(u+v)+\varphi(u-v)=2[\varphi(u)+\varphi(v)],
    \quad u,
    v\in Y.$$
Denote by $\Gamma(X)$ the set of Gaussian distributions on the group
    $X$. We note that according to this definition the degenerate distributions are Gaussian.

\bigskip

\noindent{\textbf{2. Group analogue of A.M. Kagan's theorem}}

\bigskip

\noindent Let $X$ be a second countable locally compact Abelian group. Consider a random vector    $\boldsymbol{\zeta}$= $(\zeta_1, \dots, \zeta_m)$ in values in the group $X^m$. Following A.M. Kagan \cite{Kag} we say that the distribution of the random vector $\boldsymbol{\zeta}$   belongs  to the class $\mathcal{D}_{m, k}$, $1\le k\le m$,  if the characteristic function $\hat\mu_{\boldsymbol{\zeta}}(y_1, \dots,  y_m)$ admits the following factorisation
$$
\hat\mu_{\boldsymbol{\zeta}}(y_1, \dots,  y_m)={  \mathbb{E}}[(\boldsymbol{\zeta}, (y_1, \dots,  y_m))]=\prod_{{i_1, \dots, i_k}} R_{i_1, \dots i_k}(y_{i_1}, \dots y_{i_k}), \quad y_i\in Y,
$$
where $R_{i_1, \dots i_k}(y_{i_1}, \dots, y_{i_k})$ are continuous functions such that  $R_{i_1, \dots i_k}(0, \dots, 0)=1$, and in the product  all  indexes $(i_1, \dots, i_k)$ satisfy the conditions  $1\le i_1<\dots< i_k\le m$.

Let $\xi_i$, $i=1, 2, \dots, n$,  be independent random variables with values in the group   $X$. Consider the linear forms
$L_j=\alpha_{j1}\xi_1+\cdots+\alpha_{jn}\xi_n$, $j=1, 2, \dots, m$, where the coefficients $\alpha_{ji}$ are continuous endomorphisms of the group   $X$.  It is easy to see that if these linear forms are independent, then the distribution of the random vector  $\boldsymbol{L}$= $(L_1, \dots, L_m)$  belongs to the class $\mathcal{D}_{m,  1}$.

The main result of the note is the following theorem. It is a generalisation of A.M. Kagan's result on locally compact Abelian groups.

\begin{theorem}\label{th1}.  {\it Let $X$ be a second countable locally compact Abelian group containing no subgroup topologically isomorphic to the circle group   $\mathbb{T}$.  Let $\alpha_{ji}$, $j=1, 2, \dots, m$, $i=1, 2, \dots, n$,  be continuous monomorphisms of the group  $X$.  Put $$G_i= \{(\alpha_{1i}x, \dots, \alpha_{mi}x)\in X^m: x\in X\}.$$ Assume that the following condition  hold:
\begin{equation}\label{11}
 \overline {G_i} \cap \overline {G_l}=\{0\}, \ i, l \in \{1, 2, \dots, n\}, \ i\ne l.
\end{equation}
 Let $\xi_i$, $i=1, 2, \dots, n$,  be independent random variables with values in the group   $X$ and distributions  $\mu_i$ with non vanishing characteristic functions. Consider the linear forms
$L_j=\alpha_{j1}\xi_1+\cdots+\alpha_{jn}\xi_n$, $j=1, 2, \dots, m$. If the distribution of the random vector  $\boldsymbol{L}$=$(L_1, \dots, L_m)$   belongs to the class  $\mathcal{D}_{m, m-1}$, then all $\mu_i$ are Gaussian distributions.}
\end{theorem}
For ease of reference we formulate  the
following property of adjoint homomorphisms as a lemma.
\begin{lemma}\label{lem6} {\rm (\cite[(24.41)]{Hewitt-Ross})}. {\it  Let $X_1$ and $X_2$
be locally compact Abelian groups  and  $Y_1$ and $Y_2$ be their character groups, respectively. Let $f : X_1 \mapsto X_2$ be a continuous homomorphism. Then the homomorphism $f$ satisfies $\tilde{\tilde f}=f$. The homomorphism  $\tilde f: Y_2\mapsto Y_1$ is a monomorphism if and only if the subgroup $f(X_1)$ is dense in $X_2$, and the subgroup   $\tilde f(Y_2)$ is dense in  $Y_1$ if and only if the homomorphism $f$ is a monomorphism.}
\end{lemma}

\begin{lemma}\label{lem1}. {\it  Let $X$ be a locally compact Abelian group and $Y$ be its character group.    Let $\alpha_{j}$, $j=1, 2, \dots, m$, be continuous monomorphisms of the group $X$.   Put  $a_{j}=\tilde \alpha_j$.   Let $\psi(y)$ be a continuous function on the group $Y$  satisfying the equation
\begin{equation}\label{1}
 \psi\big(a_{1}y_1+\dots +a_{m}y_m\big)= \sum_{j=1}^{m}r_j(y_1,   \dots, y_{j-1},y_{j+1},\dots, y_m), \quad y_j\in Y,
\end{equation}
where  $r_j$ are arbitrary functions. Then $\psi(y)$ is a polynomial on  $Y$.}
\end{lemma}
{\it Proof}. We use the finite difference method.  Let $h_m$ be an arbitrary element of the
group $Y$.   Substitute $y_m+h_m$  for $y_m$   in
Equation (\ref{1}). Subtracting Equation  (\ref{1})  from the resulting equation we get an equation where the right-hand side does not contain the function   $r_m$. By repeating this operation, we consistently exclude all functions $r_j$ from the right-hand side of the resulting equations.   After $m-1$ steps we get
\begin{equation}\label{2}
\Delta_{a_1h_1}\dots \Delta_{a_mh_m}\psi\big(a_{1}y_1+\dots +a_{m}y_m\big)=0, \quad y_j\in Y.
\end{equation}
By Lemma \ref{lem6}, the subgroup  ${a_{j}(Y)}$ is dense in $Y$. Since  $h_j$ are arbitrary elements of the group $Y$, it follows from (\ref{2}) that the function   $\psi(y)$ satisfies the equation
 $$
\Delta^m_{ h } \psi\big(y\big)=0, \quad   y, h\in Y,
$$
i.e. $\psi(y)$ is a polynomial. $\Box$

\begin{lemma}\label{lem2}. {\it  Let $X$ be a locally compact Abelian group and $Y$ be its character group. Let $\alpha_{pi}$, $p=1, 2, \dots, m$, $i=1, 2, \dots, n$,  be continuous monomorphisms of the group $X$. Assume that  condition    $(\ref{11})$ holds.
Put $a_{pi}=\tilde\alpha_{pi}$. Let $\psi_i(y)$ be continuous functions on the group   $Y$, satisfying the equation
\begin{equation}\label{3}
\sum_{i=1}^n\psi_i\big(a_{1i}y_1+\dots +a_{mi}y_m\big)=\sum_{j=1}^{m}s_j(y_1,   \dots, y_{j-1},y_{j+1},\dots, y_m), \quad y_j\in Y,
\end{equation}
where  $s_j$ are arbitrary functions.
Then all $\psi_i(y)$ are polynomial on $Y$.}
\end{lemma}
{\it Proof}. We use the finite difference method.  Let $h_{nj}$, $j=1, 2, \dots, m$, be   elements of the group $Y$ such that
$$
a_{1n}h_{n1}+\dots + a_{mn}h_{nm}=0.
$$
 Substitute $y_j+h_{nj}$  for $y_j$ for all $j$   in
Equation (\ref{3}) and subtract Equation  (\ref{3})  from the resulting equation.  We obtain
 \begin{equation}\label{4}
 \sum_{i=1}^{n-1}\Delta_{l_{in}}\psi_i\big(a_{1i}y_1+\dots +a_{mi}y_m\big)=\sum_{j=1}^{m}\Delta_{\bar{\bf l}_{jm}}s_j(y_1,   \dots, y_{j-1},y_{j+1},\dots, y_m), \quad y_j\in Y,
\end{equation}
where $l_{in}=a_{1i}h_{n1}+\dots + a_{mi}h_{nm}$, $i=1, 2,\dots, n-1$, ${\bf{\bar l}}_{jm}=(h_{n1},   \dots, h_{n, j-1}, h_{n, j+1}, \dots, h_{nm})$, $j=1, 2, \dots, m$. The left-hand side of Equation   (\ref{4}) no longer contains the function $\psi_n$. In the second step taking elements   $h_{n-1, j}$ in such a manner that the equality
$$
a_{1,n-1}h_{n-1,1}+\dots + a_{m,n-1}h_{n-1,m}=0
$$
holds and reasoning similarly, we exclude the function $\psi_{n-1}$ from the left-hand side of Equation   (\ref{4}).
By excluding successively the functions   $\psi_i$ from the left-hand side of Equation   (\ref{4}), after $n-1$ steps  we come to an equation of the form
\begin{equation}\label{5}
 \Delta_{l_{12}}\dots\Delta_{l_{1n}}\psi_1\big(a_{11}y_1+\dots +a_{m1}y_m\big)=\sum_{j=1}^{m}r_j(y_1,   \dots, y_{j-1},y_{j+1},\dots, y_m), \quad y_j\in Y,  \quad y_j\in Y,
\end{equation}
where
\begin{equation}\label{6}
l_{1i}=a_{11}h_{i1}+\dots + a_{m1}h_{im}, \quad i=2,   \dots, n,
\end{equation}
and $r_j$ are some functions.
It should be noted that when we got Equation (\ref{5})  we chose   elements $h_{ij}$ at every step   in such a manner that the equalities
\begin{equation}\label{7}
a_{1i}h_{i1}+\dots + a_{mi}h_{im}=0, \quad i=2,   \dots, n,
\end{equation}
were fulfilled.
Consider the subgroups $G_i=\{(\alpha_{1i}x, \dots, \alpha_{mi}x)\in X^m: x\in X\}$, $i=1, 2, \dots, n,$ in $X^m$.
Denote by $g_i$ the continuous homomorphism $g_i: Y^m\mapsto Y$ of the form
 $$
g_i(y_1, \dots, y_m)=a_{1i}y_1+\dots + a_{mi}y_m, \quad i=1, 2, \dots, n.
$$
The adjoint to $g_i$ homomorphism $\tilde g_i:X\mapsto X^m$ is of the form $\tilde g_i(x)=(\alpha_{1i}x, \dots, \alpha_{mi}x)\in G_i$.
Put   $H_i={\rm Ker}  \ g_i$.
Denote by $p_i$ the natural embedding   $p_i: H_i\mapsto Y^m$, and by $g_{li}$  the restriction of $g_l$ to $H_i$.  We have $g_{li}=g_lp_i$. Consider the annihilator $A(X^m, H_i)$. It is easy to see that     $A(X^m, H_i)=\overline {G_i}.$  The character group of the subgroup   $H_i$   is topologically isomorphic to the factor-group  $X^m/A(X^m, H_i)$, and the adjoint to   $p_i$   homomorphism  $\tilde p_i$ is a factor mapping  $\tilde p_i: X^m\mapsto X/A(X^m, H_i)$. We note that $\tilde g_{li}=\tilde p_i\tilde g_l$. Take $x\in {\rm Ker} \ \tilde g_{li}$. Then $\tilde g_l(x)=(\alpha_{1l}x, \dots, \alpha_{ml}x)\in {\rm Ker} \ \tilde p_i=A(X^m, H_i)=\overline {G_i}$. It follows from   (\ref{11}) that  $x=0$ for $i\ne l$. Thus, ${\rm Ker} \ \tilde g_{li}=\{0\}$. By Lemma \ref{lem6}, the subgroup $g_{li}(H_i)$ is dense in $Y$.

Let us return to   Equation (\ref{5}).  Since $\alpha_{j1}$,   $j=1, 2, \dots, m$,   are continuous monomorphisms of the group   $X$, by Lemma   \ref{lem1}, the function
\begin{equation}\label{41}
\Delta_{l_{12}}\dots\Delta_{l_{1n}}\psi_1\big(y\big)
\end{equation}
is a polynomial on  $Y$.
Taking into account that the subgroups $g_{li}(H_i)$, $i\ne l$, are dense in $Y$ and (\ref{6}), we can suppose that in (\ref{41}) $l_{12}=\dots=l_{1n}=h,$ where $h$ is an arbitrary element of the group $Y$, i.e. the function $\Delta_h^{n-1} \psi_1(y)$ is a polynomial on $Y$. Hence, the function $\psi_1(y)$ is also a polynomial on $Y$. For the functions  $\psi_j(y)$, $j=2, 3, \dots, n$ we reason similarly.  $\Box$

We also need group analogues of the classical Cram\'er and Marcinkiewicz theorems. We will formulate them as lemmas.

\begin{lemma}\label{lem3} {\rm (\cite{Fe1978})}. {\it Let $X$ be a second countable locally compact Abelian group containing no subgroup topologically isomorphic to the circle group   $\mathbb{T}$.  Let $\gamma\in\Gamma(X)$ and $\gamma=\gamma_1*\gamma_2$, where $\gamma_j\in {\rm M^1}(X).$ Then $\gamma_j\in {\Gamma}(X).$ }
\end{lemma}
\begin{lemma}\label{lem4} {\rm (\cite{Fe1989})}. {\it Let $X$ be a second countable locally compact Abelian group containing no subgroup topologically isomorphic to the circle group   $\mathbb{T}$. Let
$\mu\in{\rm M}^1(X)$ and the characteristic function $\hat\mu(y)$
is of the form
$$
\hat\mu(y)=\exp\{\psi(y)\}, \quad  y\in Y,
$$
where $\psi(y)$  is a continuous polynomial.
  Then $\mu\in\Gamma(X)$.}
\end{lemma}

\noindent{\it Proof of Theorem $\ref{th1}$}. Let $Y$ be a character group of the group   $X$.  Let  $\hat\mu_{\boldsymbol{L}}(y_1, \dots, y_m)$ be the characteristic function of the random vector    $\boldsymbol  L$. Put $a_{pi}=\tilde\alpha_{pi}$, $p=1, 2, \dots, m$, $i=1, 2, \dots, n$.  Taking into account that the random variables $\xi_i$ are independent,  the characteristic function $\hat\mu_{\boldsymbol{L}}(y_1, \dots, y_m)$  is of the form
$$
\hat\mu_{\bf{\boldsymbol  L}}(y_1, \dots, y_m)={  \mathbb{E}}[({\boldsymbol L},(y_1, \dots, y_m))]={  \mathbb{E}}\bigg[\prod_{i=1}^m(\alpha_{j1}\xi_1+\cdots+\alpha_{jn}\xi_n, y_j)\bigg]$$$$={  \mathbb{E}}\bigg[\prod_{i=1}^n(\xi_i, a_{1i}y_1+\dots+a_{mi}y_m)\bigg]=\prod_{i=1}^n\hat\mu_i(a_{1i}y_1+\dots+a_{mi}y_m), \quad y_j\in Y.
$$
Taking into account that the distribution of the random vector  $\boldsymbol L$  belongs to the class $\mathcal{D}_{m, m-1}$, we may write
\begin{equation}\label{9}
\prod_{i=1}^n\hat\mu_i(a_{1i}y_1+\dots+a_{mi}y_m)=\prod_{j=1}^{m} R_j(y_1,   \dots, y_{j-1},y_{j+1},\dots, y_m), \quad y_j\in Y,
\end{equation}
where   $R_j(y_1,   \dots, y_{j-1},y_{j+1},\dots, y_m)$ are continuous functions such that  $R_j(0, \dots, 0)=1$.

Put $\nu_i=\mu_i*\bar\mu_i$. Then we have $\hat\nu_i(y)=|\mu_i(y)|^2>0$,  the characteristic functions of the distributions   $\hat\nu_i(y)$ also  satisfy Equation (\ref{9}), and all the factors in the left-hand side and  right-hand side of Equation   (\ref{9}) are greater then zero. If we prove that all $\nu_i\in\Gamma(X),$ it follows from this by Lemma \ref{lem3}, that all $\mu_i\in\Gamma(X)$ too. Thus, we may assume from the   beginning that  all factors in the left-hand side and  right-hand side of  Equation   (\ref{9}) are greater then zero.

Put $\psi_i(y)=\log \hat\mu_i(y)$, $i=1, 2, \dots, n$, $s_j=\log R_j$, $j=1, 2, \dots, m$. It follows from  (\ref{9}) that the functions   $\psi_i$ and $s_j$ satisfy Equation  (\ref{3}). Then by Lemma \ref{lem2},  $\psi_i(y)$ are polynomial on $Y$. Applying Lemma \ref{lem4}, we obtain that $\mu_i\in\Gamma(X)$, $i=1, 2\dots, n$.
$\Box$

\bigskip

\noindent{\textbf{3. Corollaries and generalizations of the main theorem}}

\bigskip

\noindent Assume that   $X$ is a locally compact Abelian torsion free group. Then $X$ does not contain a subgroup topologically isomorphic to the circle group   $\mathbb{T}$, and multiplication by a nonzero integer is a continuous monomorphism of the group    $X$. Hence, Theorem   \ref{th1} implies the following statement  (compare with \cite{L1996}).
\begin{corollary}\label{co1}.  {\it Let $X$ be a second countable locally compact Abelian torsion free group.  Let $k_{ji}$, $j=1, 2, \dots, m$, $i=1, 2, \dots, n$,  be nonzero integers.  Put $G_i= \{(k_{1i}x, \dots, k_{mi}x)\in X^m: x\in X\}.$ Assume that condition   $(\ref{11})$ holds.
Let $\xi_i$, $i=1, 2, \dots, n$, be independent random variables with values in the group   $X$ and distributions  $\mu_i$ with non vanishing characteristic functions.  Consider the linear forms
$L_j=k_{j1}\xi_1+\cdots+k_{jn}\xi_n$, $j=1, 2, \dots, m$. If the distribution of the random vector  $\boldsymbol{ L}$=$(L_1, \dots, L_m)$   belongs to the class  $\mathcal{D}_{m, m-1}$, then all $\mu_i$ are Gaussian distributions.}
\end{corollary}
Suppose that in Theorem \ref{th1}  $\alpha_{ji}$, $j=1, 2, \dots, m$, $i=1, 2, \dots, n$,  are topological automorphisms of the group   $X$. Then $G_i$ are closed subgroups of the group $X^m$, and condition  (\ref{11}), as easy to see, is equivalent to the following condition:

 For any fixed $i, l\in\{1, 2, \dots, n\}$, $i\ne l$, the following equality holds:
\begin{equation}\label{12}
\bigcap_{p, k=1}^m {\rm Ker} \ (\alpha_{ki}^{-1}\alpha_{kl}-\alpha_{pi}^{-1}\alpha_{pl})=\{0\}.
\end{equation}
We note that if $m=2$, and  $L_1=\xi_1+\dots+\xi_n$ and $L_2=\alpha_1\xi_1+\dots+\alpha_n\xi_n$, where $\alpha_i$  are topological automorphisms of the group $X$, then  condition   (\ref{12}) takes  the form: $${\rm Ker} \ (\alpha_i-\alpha_l)=\{0\}$$ for all $i, l\in \{1, 2, \dots, n\}$, $i\ne l$.

Let $X$  be a locally compact Abelian divisible torsion free group.  We recall that an Abelian group $G$ is divisible if, for every natural $n$ and every $g\in G$, there exists $x\in G$ such that $nx=g$.
 An important example of a locally compact Abelian divisible torsion free group is a group of the form   $\mathbb{R}^n\times{\Sigma^{\mathfrak{n}}_{\text{\boldmath $a$}}}$, where  $\Sigma_{\text{\boldmath $a$}}$ is an \text{\boldmath $a$}-adic solenoid  (a compact Abelian group such that its character group is topologically isomorphic to the group of rational numbers $\mathbb{Q}$ considering   the discrete topology),    ${\mathfrak{n}}$ is a cardinal number. Any connected locally compact Abelian divisible torsion free group is topologically isomorphic
 to a group of this type. For the structure theorem for an arbitrary locally compact Abelian divisible
 torsion free group see  \cite[(25.33)]{Hewitt-Ross}. Obviously, for these groups $X$ multiplication by any nonzero integer, and hence, multiplication by any nonzero rational number, is a topological automorphism.  Therefore, we may consider  linear forms of independent random variables taking values in   $X$ with
  rational coefficients. Corollary \ref{co1} for such groups  can be significantly improved. Namely, condition  (\ref{11}) or the equivalent  condition (\ref{12}),   can be omitted.
\begin{theorem}\label{th3}.  {\it Let $X$  be a second countable locally compact Abelian divisible torsion free group.  Let $\alpha_{ji}$, $j=1, 2, \dots, m$, $i=1, 2, \dots, n$,   be nonzero rational numbers. Let $\xi_i$, $i=1, 2, \dots, n$,  be independent random variables with values in the group   $X$ and distributions  $\mu_i$ with non vanishing characteristic functions. Consider the linear forms
$L_j=\alpha_{j1}\xi_1+\cdots+\alpha_{jn}\xi_n$, $j=1, 2, \dots, m$. If the distribution of the random vector  $\boldsymbol{L}$=$(L_1, \dots, L_m)$   belongs to the class  $\mathcal{D}_{m, m-1}$, then all $\mu_i$ are Gaussian distributions.}
\end{theorem}
{\it  Proof}.  It is easy to see that the character group of a locally compact Abelian divisible torsion free group is also a locally compact Abelian divisible torsion free group. The proof of Theorem \ref{th1} is based on  Lemma  \ref{lem2} and group analogues of theorems by Cram\'er and Marcinkiewicz (Lemmas \ref{lem3} and \ref{lem4}).  Since $X$ does not contain a subgroup topologically isomorphic to the circle group   $\mathbb{T}$, Lemmas \ref{lem3} and \ref{lem4} hold for the group $X$.   We retain the notation used in the proof of the Theorem \ref{th1}.    Reasoning as in the proof of Theorem    \ref{th1}, we come to   Equation  (\ref{3}).

Consider integer-valued vectors ${\boldsymbol{ \alpha}}_i=(\alpha_{1i}, \dots, \alpha_{mi})$, $i=1, 2, \dots, n$. Since multiplication by any nonzero rational number is a topological automorphism of the groups   $X$ and $Y$, it is easy to see that condition  (\ref{11}) or equivalent condition  (\ref{12}), are not fulfilled if and only if the vectors   ${\boldsymbol{ \alpha}}_i$ and  ${\boldsymbol{ \alpha}}_l$ are   collinear for some $i\ne l$. We assume for simplicity that we have the only subset of collinear vectors, and they are the vectors      ${\boldsymbol{ \alpha}}_q, \dots, {\boldsymbol{ \alpha}}_{n}$, and the vectors ${\boldsymbol{ \alpha}}_{1}, \dots, {\boldsymbol{ \alpha}}_{q}$  are not collinear. The general case may be considered similarly. Let   ${\boldsymbol{ \alpha}}_{i}=c_i{\boldsymbol{\ \alpha}}_{q}$, $i=q+1, q+2, \dots, n$, where $c_i$ are nonzero rational numbers. Put
$$\varphi_i(y)=\psi_i(y),  \ i=1, 2, \dots, q-1, \quad \varphi_q(y)=\psi_q(y)+\psi_q(c_{q+1}y)+\dots+ \psi_q(c_{n}y).$$
Equation (\ref{3}) with this notation takes the form
\begin{equation}\label{31}
\sum_{i=1}^q\varphi_i\big(\alpha_{1i}y_1+\dots +\alpha_{mi}y_m\big)=\sum_{j=1}^{m}s_j(y_1,   \dots, y_{j-1},y_{j+1},\dots, y_m), \quad y_j\in Y.
\end{equation}
The non-collinearity of the vectors   ${\boldsymbol{ \alpha}}_{i}$, $i=1, 2, \dots, q$, implies that the topological automorphisms   $\alpha_{pi}$, $p=1, 2, \dots, m$, $i=1, 2, \dots, q$,  already satisfy   condition    (\ref{11})  or equivalent condition    (\ref{12}). Now,   the conditions of    Lemma  \ref{lem2} are satisfied. By Lemma \ref{lem2}, the functions $\varphi_i(y)$, $i=1, 2, \dots, q$,  are polynomial on $Y$. Applying Lemma \ref{lem4}, we obtain that $\mu_i\in\Gamma(X)$, $i=1, 2,\dots, q-1$. Moreover, $\exp\{\varphi_q(y)\}$ is the characteristic function of a Gaussian distribution. Taking into account that  multiplication by   $c_i$, $i=q+1, q+2, \dots, n$, are topological automorphisms of the group   $Y$, and applying Lemma \ref{lem3}, we get that $\mu_i\in\Gamma(X)$, $i=q, q+1,\dots, n$.

It should be noted that if in Theorem \ref{th3}  $m=2$, i.e. we have two independent linear forms   $L_1$ and $L_2$ of $n$ independent random variables $\xi_i$, and the coefficients of the linear forms   $\alpha_{ji}$, $j=1, 2, \dots, m$, $i=1, 2, \dots, n$, are nonzero integers, then, as it was proved in   \cite{Fe1990}, the statement of Theorem \ref{th3} remains true if we omit the condition that $X$ is a divisible group.
$\Box$

In the article \cite{KS} A. M. Kagan and G. J. Sz\'ekely introduced
a notion of $Q$-independence and proved, in particular,
that the     Skitovich--Darmois theorem
holds true if instead of independence  $Q$-independence is considered.
  Then in the article \cite{F2017} a notion of $Q$-independence for random variables with values in a locally compact Abelian group  was introduced.  In \cite{F2017} it was proved, in particular, that a group analogue of the Skitovich--Darmois theorem (see  \cite{Fe-SD-2003}) remains true if instead of independence  $Q$-independence is considered.

Let $\xi_1, \dots, \xi_n$ be random variables with values in the group
 $X$.
We say that the random variables
  $\xi_1,  \dots, \xi_n$ are
$Q$-independent  if   the characteristic function of the random vector   $\boldsymbol{\xi}$=$(\xi_1, \dots, \xi_n)$ is represented in the form
\begin{equation}\label{17}
\hat\mu_{\boldsymbol{\xi}}(y_1, \dots, y_n)={  \mathbb{E}}[(\xi_1, y_1)\cdots(\xi_n, y_n)]=$$$$=\bigg(\prod_{j=1}^n\hat\mu_{\xi_j}(y_j)\bigg)\exp\{q(y_1,
\dots, y_n)\}, \quad y_j\in Y,
\end{equation}
where $q(y_1, \dots, y_n)$ is a continuous polynomial on the group   $Y^n$.
We   also suppose that   $q(0, \dots, 0)=0$.

We will prove now that, exactly as in the case of the Skitovich--Darmois theorem, Theorem    \ref{th1} holds true if we change the condition of independence of the random variables $\xi_i$ for $Q$-independence. The following statement holds.

\begin{theorem}\label{th4}.  {\it Let $X$ be a second countable locally compact Abelian group containing no subgroup topologically isomorphic to the circle group   $\mathbb{T}$.  Let $\alpha_{ji}$, $j=1, 2, \dots, m$, $i=1, 2, \dots, n$,  be continuous monomorphisms of the group  $X$.  Put $G_i= \{(\alpha_{1i}x, \dots, \alpha_{mi}x)\in X^m: x\in X\}.$ Assume that   condition  $(\ref{11})$ holds.
 Let $\xi_i$, $i=1, 2, \dots, n$,  be $Q$-independent random variables with values in the group   $X$ and distributions  $\mu_i$ with non vanishing characteristic functions. Consider the linear forms
$L_j=\alpha_{j1}\xi_1+\cdots+\alpha_{jn}\xi_n$, $j=1, 2, \dots, m$. If the distribution of the random vector  $\boldsymbol{L}$=$(L_1, \dots, L_m)$   belongs to the class  $\mathcal{D}_{m, m-1}$, then all $\mu_i$ are Gaussian distributions.}
\end{theorem}

\noindent{\it Proof}. We retain the notation used in the proof of  Theorem \ref{th1}, and reason as in the proof of Theorem    \ref{th1}.   Taking into account  the $Q$-independence of random variables $\xi_i$, the characteristic function of the random vector  $\boldsymbol{L}$ is of the form
$$\hat\mu_{\boldsymbol{L}}(y_1, \dots, y_m)=\bigg(\prod_{i=1}^n\hat\mu_i(a_{1i}y_1+\dots+a_{mi}y_m)\bigg)
$$$$\times\exp\{q(a_{11}y_1+\dots+a_{m1}y_m,
\dots, a_{1n}y_1+\dots+a_{mn}y_m)\}), \quad y_j\in Y,
$$
where $q(y_1, \dots, y_n)$ is a continuous polynomial on the group   $Y^n$, and instead of Equation  (\ref{9}) we get the equation
\begin{equation}\label{51}
\prod_{i=1}^n\hat\mu_i(a_{1i}y_1+\dots+a_{mi}y_m)$$$$ =\bigg(\prod_{j=1}^{m} R_j(y_1,   \dots, y_{j-1},y_{j+1},\dots, y_m)\bigg)\exp\{p(y_1, \dots, y_m)\}), \ \ y_j\in Y,
\end{equation}
where $p(y_1, \dots, y_m)$ is a continuous polynomial on the group   $Y^m$. Following the proof of Theorem \ref{th1}, we want to conclude from (\ref{51}) that all distributions $\mu_i$ are Gaussian. To do this we need to know that the statements of Lemma  \ref{lem1} and \ref{lem2} remain true in the case when we add a polynomial to the right-hand side of Equations   (\ref{1})  and (\ref{3}). It is easy to see that it is true. The proofs of the lemmas are almost the same. The only difference is that the degrees of polynomials   $\psi(y)$  and $\psi_i(y)$ are changed. $\Box$

Theorem \ref{th4} and the reasoning used in the proof of Theorem \ref{th3} imply  the following analogue of Theorem \ref{th3} for  Q-independent random variables.
\begin{theorem}\label{th5}.  {\it Let $X$  be a second countable locally compact Abelian divisible torsion free group.  Let $\alpha_{ji}$, $j=1, 2, \dots, m$, $i=1, 2, \dots, n$,   be nonzero rational numbers. Let $\xi_i$, $i=1, 2, \dots, n$,  be Q-independent random variables with values in the group   $X$ and distributions  $\mu_i$ with non vanishing characteristic functions. Consider the linear forms
$L_j=\alpha_{j1}\xi_1+\cdots+\alpha_{jn}\xi_n$, $j=1, 2, \dots, m$. If the distribution of the random vector  $\boldsymbol{L}$=$(L_1, \dots, L_m)$   belongs to the class  $\mathcal{D}_{m, m-1}$, then all $\mu_i$ are Gaussian distributions.}
\end{theorem}

\begin{remark}\label{re2}.
{\rm  Assume that in Theorem \ref{th1}  $m=2$, i.e. we have two independent forms $L_1$ and $L_2$ of $n$ independent random variables   $\xi_i$. Suppose that  $\alpha_{ji}$, $j=1, 2$, $i=1, 2, \dots, n$,  are topological automorphisms of the group   $X$. Then, as it was proved in   \cite{Fe-SD-2003}, all $\mu_i$ are Gaussian distributions. In so doing we  do not suppose that   condition    (\ref{11})   or equivalent condition    (\ref{12}) is valid.

Considering this we will formulate the following question. Does Theorem
\ref{th1} hold true if we omit
condition     (\ref{11})   or  (\ref{12}) in the case when the coefficients $\alpha_{ji}$   are topological automorphisms of the group  $X$?
}
\end{remark}
\begin{remark}\label{re3}.
{\rm   Discuss the case  when a locally compact Abelian group $X$ contains a subgroup topologically isomorphic to the circle group $\mathbb{T}$. The simplest example of such a group  is $X=\mathbb{T}$. Then, as it was proved in  \cite{Fe-Be-1986}, there exist independent non-Gaussian random variables   $\xi_1$ and $\xi_1$ with values in the group  $X$ with nonvanishing characteristic functions such that the linear forms    $L_1=\xi_1+\xi_2$ and $L_2=\xi_1-\xi_2$ are independent, i.e. the distribution of the random vector  ${\boldsymbol L}=(L_1,   L_2)$   belongs to the class $\mathcal{D}_{2,  1}$. Condition  $(\ref{12})$ in this case takes  the form: multiplication by 2 is a topological automorphism of the group $X$. Obviously, it is not true.

Consider a more complicated example. Let $X=\mathbb{T}^2$. The main theorem proved in   \cite{MiFe4} implies the complete description of topological automorphisms  $\alpha$ of the group $X$ which possess the following property: there exist independent non-Gaussian random variables   $\xi_1$ and $\xi_1$ with values in the group  $X$ with nonvanishing characteristic functions such that the linear forms   $L_1=\xi_1+\xi_2$ and $L_2=\xi_1+\alpha\xi_2$ are independent, i.e. the distribution of the random vector  ${\boldsymbol L}=(L_1,   L_2)$   belongs to the class $\mathcal{D}_{2,  1}$.
Condition  $(\ref{12})$ in this case takes  the form: $I-\alpha$ is a topological automorphism of the group $X$. This condition for such  $\alpha$ is not true. As follows from the results of   \cite{MyFe2011}, the same is also true for the group $X=\mathbb{R}\times\mathbb{T}$.

In connection with the above, the following question arises. Let $X$ be an arbitrary second countable locally compact Abelian group, i.e. generally speaking, $X$ can contain a subgroup topologically isomorphic to the circle group   $\mathbb{T}$.  Assume that condition  (\ref{11})   or  (\ref{12}) in the case when coefficients   $\alpha_{ji}$   are topological automorphisms of the group  $X$, is valid. Does Theorem \ref{th1} hold true in this case? In other words, does  Theorem \ref{th1} hold true  if we omit the condition: the group $X$ contains no subgroup topologically isomorphic to the circle group $\mathbb{T}$. It is interesting to remark that if $m=2$, i.e. we have two independent linear forms, the number of independent random variables $n=2$,   and the coefficients    $\alpha_{ji}$ are topological automorphisms of the group $X$,   as it was proved in   \cite{Fe2009}  the answer too this question is positive.
}
\end{remark}


 \newpage

\noindent B. Verkin Institute for Low Temperature Physics and Engineering \\ of the National Academy
of Sciences of Ukraine\\ Nauky Ave. 47,\\
Kharkiv, 103,\\ 61103, Ukraine

\medskip

\noindent e-mail: feldman@ilt.kharkov.ua



\newpage


\begin{thebibliography}{99}

\bibitem{Fe1978} G.M. Feldman: Gaussian distributions on locally compact abelian groups.   Theory Probab. Appl. \textbf{23} (3),  529--542  (1979)

\bibitem{Fe-Be-1986} G.M. Feldman: Bernstein Gaussian distributions on
groups.  Theory Probab. Appl.  \textbf{31} (1),   40--49 (1987)

\bibitem{Fe1989} G.M. Feldman: Marcinkiewicz and Lukacs theorems on
abelian groups.  Theory Probab. Appl.  \textbf{34} (2),   290–-297 (1990)

\bibitem{Fe1990} G.M. Feldman: Characterization of the
    Gaussian distribution on groups by the independence of
linear statistics. Siberian Math. J. \textbf{31} (2),   336--345 (1990)

\bibitem{Fe1996} G.M. Feldman: On the Skitovich-Darmois theorem on
    compact groups. Theory Probab. Appl.     \textbf{ 41} (4),    768--773 (1997)

\bibitem{Fe1997} G.M. Feldman: The Skitovich-Darmois theorem for discrete
periodic Abelian groups.  Theory Probab. Appl.  \textbf{42} (4),    611--617 (1998)

\bibitem{Fe-SD-2003} G.M. Feldman: A characterization of the Gaussian
distribution on Abelian groups.  Probab. Theory Relat. Fields.
 \textbf{126} (1),   91--102 (2003)

 \bibitem{Fe1} G.M. Feldman: Functional equations and characterization problems on locally compact Abelian groups. EMS Tracts in Mathematics \textbf{5}, Zurich. European Mathematical Society (EMS), 2008

\bibitem{Fe2009} G.M. Feldman: On a theorem of K. Schmidt. Bull. Lond. Math. Soc. \textbf{41} (1),   103--108 (2009)

\bibitem{F2017}  G.M. Feldman:  Characterization theorems for $Q$-independent random variables with values in a locally compact Abelian group. Aequationes Mathematicae. \textbf{91},  949--967 (2017)

\bibitem{FG2005} G.M. Feldman, P. Graczyk: The Skitovich-Darmois theorem for locally compact Abelian groups. J. of the Australian Mathematical Society.  \textbf{88} (3),  339--352 (2010)


\bibitem{MyFe2011} G.M. Feldman, M.V. Myronyuk: Independent linear forms on connected Abelian groups. Math. Nachr.  \textbf{284} (2–3),  255--265 (2011)

    \bibitem{Hewitt-Ross}  E. Hewitt,  K.A. Ross: Abstract
   Harmonic Analysis  {\bf 1},  Springer-Verlag, Berlin,
    Gottingen, Heildelberg, 1963

\bibitem{Kag} A.M. Kagan: New Classes of dependent random variables and a generalization of the Darmois--Skitovich theorem to several forms.
Theory Probab. Appl. \textbf{33} (2),   286--295 (1989)


\bibitem{Kag-Lin-Rao}   A.M. Kagan, Yu. V. Linnik,   C.R. Rao:
{Characterization problems in mathematical statistics}, Wiley Series in Probability and Mathematical
Statistics, John Wiley $\&$ Sons, New York-London-Sydney, 1973

\bibitem{KS} {  A.M. Kagan, G.J. Sz\'ekely:} An analytic generalization of independence and identcal distributiveness. Statistics and Probability Letters  \textbf{110},   244--248 (2016)

    \bibitem{L1996} S.V. Lisyanoi:  On a group analogue of a theorem of     A. M. Kagan.  Theory Probab. Appl.  \textbf{40} (1),    165–-167 (1995)

\bibitem{MiFe4}  M.V. Myronyuk, G.M. Feldman: Independent linear statistics
on the two-dimensional torus. Theory Probab. Appl.  \textbf{52} (1),  78--92 (2008)


\bibitem{Parthasarathy} K.R. Parthasarathy: Probability measures
on metric spaces.   New York and London: Academic Press, 1967










\end{thebibliography}
\end{document}